\documentclass[11pt]{amsart}

\usepackage{amscd}
\usepackage{amsmath, amssymb,color, stmaryrd}
\usepackage{amsfonts}
\usepackage{comment}

\newcommand{\de}{\partial}

\newcommand{\ddbar}{\sqrt{-1}\partial\overline{\partial}}
\newcommand{\p}{\partial}
\newcommand{\pbar}{\ov{\partial}}

\newcommand{\norm}[1]{\left\lVert#1\right\rVert}

\newcommand{\ov}[1]{\overline{#1}}
\newcommand{\ul}[1]{\underline{#1}}

\newcommand{\vp}{\varphi}
\newcommand{\e}{\varepsilon}

\newcommand{\Sing}{\mathrm{Sing}}

\newcommand{\mSH}{m\mathrm{SH}}
\newcommand{\SH}{\mathrm{SH}}

\newcommand{\Hm}{\mathrm{H}^m}

\newcommand{\Supp}{\mathrm{Supp}}

\newcommand{\R}{\mathbb{R}}
\newcommand{\C}{\mathbb{C}}
\newcommand{\CP}{\mathbb{CP}}

\renewcommand{\leq}{\leqslant}
\renewcommand{\geq}{\geqslant}

\newcommand{\be}{\begin{equation}}
\newcommand{\ee}{\end{equation}}

\begin{document}
\newtheorem{claim}{Claim}
\newtheorem{theorem}{Theorem}[section]
\newtheorem{lemma}[theorem]{Lemma}
\newtheorem{corollary}[theorem]{Corollary}
\newtheorem{proposition}[theorem]{Proposition}
\newtheorem{question}{question}[section]
\theoremstyle{definition}
\newtheorem{definition}[theorem]{Definition}
\newtheorem{remark}[theorem]{Remark}

\numberwithin{equation}{section}

\title[Lelong numbers of $m$-subharmonic functions]{Lelong Numbers of $m$-Subharmonic Functions Along Submanifolds}

\begin{abstract}
We study the possible singularities of an $m$-subharmonic function $\varphi$ along a complex submanifold $V$ of a compact K\"ahler manifold, finding a maximal rate of growth for $\varphi$ which depends only on $m$ and $k$, the codimension of $V$. When $k < m$, we show that $\varphi$ has at worst log poles along $V$, and that the strength of these poles is moveover constant along $V$. This can be thought of as an analogue of Siu's theorem.
\end{abstract}

\author[J. Chu]{Jianchun Chu}
\address[J. Chu]{School of Mathematical Sciences, Peking University, Yiheyuan Road 5, Beijing, P.R.China, 100871}
\email{jianchunchu@math.pku.edu.cn}

\author[N. McCleerey]{Nicholas McCleerey}
\address[N. McCleerey]{Department of Mathematics, Univeristy of Michigan, Ann Arbor, 530 Church St, Ann Arbor, MI 48109}
\email{njmc@umich.edu}


\maketitle


\section{Introduction}


Let $(X^{n},\omega)$ be a closed K\"ahler manifold, and let $\mSH(X, \omega)$ be the space of $(m, \omega)$-sh functions on $X$, $1 \leq m\leq n$. $m$-sh functions were introduced by B\l ocki \cite{Bl05} as the natural space of weak solutions to the complex $m$-Hessian equation on domains in $\C^n$, and they naturally interpolate between plurisubharmonic (psh) and subharmonic (sh) functions as $m$ varies.

As such, $m$-sh functions share many similarities with psh functions. Both support a robust potential theory which has been the subject of much interest in recent years (see e.g. \cite{BT76, BT82, Bl05, ACCH09, DK18, LN19, DDL21}). One major difference however is the lack of a strong relation between $m$-sh functions and analytic geometry, something which can be seen, for instance, in our inability to solve the $\pbar$-problem with generic $m$-sh weights. This geometric connection is utilized to great effect in the study of psh functions, and is largely responsible for our understanding what singularities psh function can look like e.g. Demailly approximation \cite{Dem92} and Siu's theorem \cite{Siu74}.

Our goal in this paper is to develop a better understanding of the singularities of $m$-sh functions, proceeding in a more ad-hoc manner. A natural starting place is to study their behaviour along complex submanifolds:

\begin{theorem}\label{psi_V_intro}
Let $V\subset X$ be a closed, complex submanifold of codimension $k$ and $\Omega$ a sufficiently small $m$-hyperconvex neighborhood of $V$.

For each $m \leq k$, there exists an $(m,\omega)$-sh function $\psi_V$ on $X$, locally bounded and maximal on $\Omega\setminus V$, with the same singularity type as:
\[
G_m(r) := \begin{cases} \log r &\text{ when }m=k\\ -r^{-2\left(\frac{k}{m}-1\right)} & \text{ when }m < k,\end{cases}
\]
where here $r := \mathrm{dist}_\omega(\cdot, V)$ is the $\omega$-distance to $V$.
\end{theorem}

With the function $\psi_V$ in hand, we can study generalized Lelong numbers, in the sense of Demailly \cite{Dem93}, and relative types, in the sense of Rashkovskii \cite{Rash06}, of arbitrary $m$-sh functions along $V$.

\begin{corollary}\label{Lelong_number_intro}
Let $V$, $\psi_V$ and $r$ be as in Theorem \ref{psi_V_intro}. Then for any $\vp\in\mSH(X, \omega)$, the limit:
\[
\nu(\vp, \psi_V) := \lim_{s\rightarrow 0} \frac{C_{k,m}}{s^{2k - 2k/m}} \int_{\{r < s\}} \ddbar \vp \wedge (\ddbar r^2)^{m-1}\wedge\omega^{n-m}
\]
exists, and is both finite and non-negative. Here $C_{k,m}$ is a constant depending only on $k$ and $m$.
\end{corollary}

\begin{corollary}\label{relative_type_intro}
Let $V$ and $\psi_V$ be as in Theorem \ref{psi_V_intro}. Let $W_{s} := \{\psi_V < s\}$. Then for any $\vp\in\mSH(X, \omega)$, the limit:
\[
\sigma(\vp, \psi_V) := \lim_{s\rightarrow-\infty} \frac{\max_{W_s} \vp}{s},
\]
exists, and is both finite and non-negative.
\end{corollary}

Note that, if $m > k$, then we can still apply Corollaries \ref{Lelong_number_intro} and \ref{relative_type_intro} by just considering $\vp \in \mSH(X,\omega)\subset k\SH(X,\omega)$.

Corollaries \ref{Lelong_number_intro} and \ref{relative_type_intro} are the first such results about the singularities of $m$-sh functions along submanifolds of positive dimension -- for more information about Lelong numbers of $m$-sh functions at points, we refer the reader to \cite{DK18}, and also \cite{HL_part1, HL_part2, Chu21}. As one sees from Theorem \ref{psi_V_intro} and Corollary \ref{Lelong_number_intro}, the natural scaling changes with the codimension $k$ when $k \geq m$, something quite different from the psh case.

The scaling we obtain agrees with what one might naively expect if they were to assume that the restriction of $\vp$ to any $k$-dimensional submanifold transverse to $V$ were still $m$-sh. This assumption is very much not true; $\vp$ will generally lose positivity when restricted to a submanifold, since $\ddbar \vp$ can have negative eigenvalues. That we still recover this optimal scaling may seem surprising in light of this.

This loss of positivity also means that we cannot use a slicing argument to reduce our situation to the case of studying the singularity of $\vp$ at the origin, as one does in the psh case. Instead, we show Theorem \ref{psi_V_intro} by constructing smooth sub- and super-solutions (essentially) to the equation $H^m(\vp) = 0$ on $\Omega\setminus V$ by taking small perturbations of $G_m(r)$, using the computations of Tam-Yu \cite{TY12} for the Hessian of $r$. The function $\psi_V$ is then realized as a singularity type envelope; the sub- and super-solutions guarantee that $\psi_V$ has the correct singularity type.

We show that our construction of sub-solutions can be modified in Proposition \ref{localization} to produce $m$-sh functions whose behaviour near $V$ is roughly like:
\[
-\theta G_m(r)
\]
for any smooth $\theta\geq 0$ (see Proposition \ref{localization} for a precise statement). This behaviour is predicted by \AA hag-Cegrell-Czy\.{z}-Hi\d{\^{e}}p's \cite{ACCH09}, and Hung-Phu's \cite{HP17}, solutions to certain highly degenerate complex Monge-Amp\`ere (resp. complex Hessian) equations. Their solutions are constructed as a sequence of envelopes, and the resulting behaviour of the solutions near $V$ is not easily seen. Our examples provide constraints on that behaviour can be (which can be quite strong in fact -- see Remark \ref{ACCH_remark}), but without corresponding super-solutions, we cannot deduce the exact singularity types of the solutions.

In the context of this paper, the examples in Proposition \ref{localization} are interesting in that they show that $\sigma(\vp, \psi_V) \not= \nu(\vp, \psi_V)$ in general, unlike the more classical case when the weight is singular only at a point. Indeed, these examples show that it is possible to have $\sigma(\vp, \psi_V) = 0$ but $\nu(\vp, \psi_V) > 0$ -- the reverse is impossible (see Corollary \ref{Lelong_Relative_Compare}).

We conclude that we cannot expect constraints which are significantly stronger than Corollaries \ref{Lelong_number_intro} and \ref{relative_type_intro} when $V$ has codimension $\geq m$. When the codimension is less than $m$ however, our next result finds that the situation is much more similar to the psh case:

\begin{theorem}\label{Siu_intro}
Suppose that $V\subset X$ is a complex submanifold of codimension $k$, and $\vp\in\mSH(X,\omega)$. Suppose that $k < m$, and let $\psi_V\in k\mathrm{SH}(X,\omega)$ be the function constructed in Theorem \ref{psi_V_intro}. Then the function:
\[
L_{\psi_V}(\vp)(z) := \liminf_{\substack{z'\rightarrow z\\ z'\not\in V}} \frac{\vp(z')}{\psi_V(z')}\quad z\in V,
\]
is constant along $V$.
\end{theorem}

Theorem \ref{Siu_intro} can be seen as a version of Siu's theorem in our context. Note in particular that implies a certain propagation of singularities along $V$, since the function $L_{\psi_V}(\vp)$ is defined pointwise.

Since the analyticity of $V$ is already assumed, the proof of Theorem \ref{Siu_intro} is pleasantly elementary, and boils down to essentially the weak Harnack inequality and Corollary \ref{relative_type_intro}.

\medskip

Taken together, our results indicate a stark dichotomy in the behaviour of $m$-sh functions, depending on the size of their singular set. The situation where the singular set is small appears to be easier, since it can be ``selected" by assuming that $\vp$ has a well-defined complex Hessian measure -- specifically, one sees that, if we assume that $\vp\in \mathcal{E}^m$, the $m$-subharmonic version of Cegrell's class, then $\nu(\vp, \psi_V) = 0$ for all $k \leq m$. The small codimension case is likely to me more interesting from a geometric point of view however, and it is unclear if there is a similar potential theoretic assumption one can make to limit to this situation. Finding conditions which work well in both cases appears challenging.

We conclude this introduction with an outline of the rest of the paper. In Section \ref{section_background} we recall some background results, which should be more-or-less standard to experts. Since we will need to utilize the lack of boundary of $V$ at several points, it is important that we do not restrict ourselves to domains in $\C^n$, and work rather on abstract $m$-hyperconvex manifolds. In Section \ref{section_constructions}, we prove Theorem \ref{psi_V_intro} (Theorem \ref{maximal_construction}), as well as Corollary \ref{Lelong_number_intro} (Proposition \ref{integral_Lelong_Laplacian}). The examples in Proposition \ref{localization} are constructed in Subsection \ref{subsection_localization}. Theorem \ref{Siu_intro} is finally shown in Section \ref{section_Siu} (Theorem \ref{Siu}).

\medskip

{\bf Acknowledgements:} We would like to thank M. Jonsson for his interest in the present work. J. Chu was partially supported by Fundamental Research Funds for the Central Universities (No. 7100603592 and No. 7100603624).


\section{Background and Notation}\label{section_background}


Throughout, we shall assume that $(X^n, \omega)$ is a closed K\"ahler manifold of complex dimension $n$, with K\"ahler form $\omega$. $(\Omega^n, \omega)$ will always denote a compact K\"ahler manifold with boundary. We will often assume that $\Omega\subset X$. We always assume that $V^{n-k}\subset \Omega$ is a compact submanifold without boundary and (complex) codimension $k$ (so that $V$ has complex dimension $n - k$).

\begin{definition}
Suppose that $(\Omega^n, \omega)$ is a K\"ahler manifold, with boundary, and let $m$ be an integer between $1$ and $n$. We say a smooth $(1,1)$-form $\alpha$ is $m$-subharmonic (or $m$-sh) if
\[
\alpha^k\wedge\omega^{n-k} \geq 0\text{ for each } k\in\{1, \ldots, m\}.
\]
We denote the set of all $m$-sh forms by $\Gamma^m(\Omega)$, or $\Gamma_\omega^m(\Omega)$ if the metric needs to be specified.

In \cite{Bl05}, B\l ocki defines a real $(1, 1)$-current $T$ to be $m$-sh if
\[
T\wedge\alpha_{2}\wedge\ldots\wedge\alpha_m\wedge\omega^{n-m} \geq 0\text{ for all }\alpha_{2},\ldots, \alpha_m\in\Gamma^m(\Omega).
\]
By G\aa rding's inequality \cite{Gard59}, the above definitions are consistent. Additionally, we shall say that $T$ is strictly $m$-sh if $T - \delta\omega$ is $m$-sh for some $\delta > 0$.

We say a function $\rho$ on $\Omega$ is $m$-sh if the $(1,1)$-current $dd^c\rho$ is $m$-sh, and write $\mSH(\Omega)$ for the space of $m$-sh functions on $\Omega$ (or $\mSH_\omega(\Omega)$, if the metric needs to be specified). We will also say that $\rho$ is strictly $m$-sh if $i\p\pbar\rho$ is.

Given a closed, real $(1,1)$-form $\theta$, we say that a function $u$ is $(m,\theta)$-sh if $\theta_u := \theta + dd^c u$ is $m$-sh, and write $\mSH(X, \theta)$ for the set of all $(m, \theta)$-sh functions on $X$.

Classically, the case which has been most studied is when $\theta = \omega$; however, it is easy to check that essentially all standard results hold when $\theta$ is any closed, strictly $m$-sh form.
\end{definition}


\subsection{$m$-hyperconvex manifolds}


We have the following $m$-sh analogue of hyperconvex manifolds.

\begin{definition}\label{hyperconvex}
We say a compact K\"ahler manifold $(\Omega, \omega)$ is $m$-hyperconvex if there exists a strictly $m$-sh exhaustion function on $\Omega^\circ$ which is smooth up to the boundary i.e. a strictly $m$-sh function $0 > \rho \in C^\infty(\ov{\Omega})$, such that $\Omega_c := \{\rho < c\} \Subset \Omega$ for each $c < 0$.
\end{definition}

It is easy to see that $m$-hyperconvex manifolds can contain non-trivial closed subvarieties. This basic observation is key to our setting, so it is important that we do not restrict our definition to domains in $\C^n$, which has been the primary case considered in much of the literature. Many of those previous results do not heavily utilize this assumption however, and can be shown to hold for arbitrary hyperconvex domains with only minor changes to their proofs.

One difference between the two settings which is worth pointing out is that it is crucial for $\rho$ to be strictly $m$-subharmonic in Definition \ref{hyperconvex}. Without this, the complex Hessian operator may fail to be well-defined; for example, consider a neighborhood $U$ of a smooth, ample divisor $D \subset \CP^n$. If we let $s$ be a holomorphic section of $\mathcal{O}_{\CP^n}(D)$ with $\{s = 0\} = D$, and $h$ a positive metric on $\mathcal{O}_{\CP^n}(D)$, then it is well known that the product $(\ddbar \log |s|_h)^2\wedge\omega^{n-m}$ is not well defined, even though there exists a smooth, psh exhaustion function on $U$. For positive results in this direction, see the recent papers \cite{And05, AW14, ABW19, Bl19, AWNW21}, which define product currents which retain some of the singular nature of $\log |s|_h$.\\

Since we repeatedly make use of the fact that codimension $m$-subvarieties admit $m$-hyperconvex neighborhoods, we record this fact here:

\begin{proposition}\label{mhc_nbrhd}
Suppose that $(X^n, \omega)$ is a compact K\"ahler manifold and that $V\subset X$ is a compact submanifold of codimension $m$. Then there exists an $m$-hyperconvex neighborhood $U$ of $V$ with smooth boundary.
\end{proposition}

\begin{proof}
Let $r(z) := \mathrm{dist}_\omega(z, V)$ be the Riemannian distance function to $V$. When $r$ is sufficiently small, $r^{2}$ is smooth. By \cite{TY12},
\[
\ddbar r^2 = \begin{pmatrix}I_k & 0\\ 0 & 0 \end{pmatrix} + o(1),
\]
where $o(1)$ denotes a term satisfying $\lim_{r\rightarrow0}o(1)=0$. It is clear that the leading term belongs to $\Gamma^{m}$. Then we can take any sublevel set of $r^{2}$ to be the desired $m$-hyperconvex neighborhood of $V$.
\end{proof}


\subsection{The Complex Hessian Operator}


For later use, we recall some basic facts about the complex Hessian operator on $m$-hyperconvex domains -- see e.g. \cite{WanWang16}, which draws heavily from the work of Demailly \cite{Dem85, Dem87, Dem_Potential, Dem93, Dem_CADG_book} in the psh case. The above papers phrase their results in terms of a weight function $\psi$, and they required the unbounded locus of $\psi$ to be a discrete set. However, it is easy to see that this assumption is superfluous for the results we will need below -- all that is needed is for the unbounded locus of $\psi$ to be compact and that it not intersect the boundary of $\Omega$, as in \cite{Dem_CADG_book}.\\

Throughout, we will be assuming that $(\Omega^n, \omega)$ is an $m$-hyperconvex manifold (which recall we assume to be compact with boundary). We start by defining the classes of weights which we will study:
\begin{definition}
We say that $\psi$ is a {\bf weight} on $\Omega$ if $\psi$ satisfies the following:
\begin{enumerate}
\item $\psi$ is bounded above.
\item $\Sing(\psi) := \{\psi = -\infty\}$ is closed.
\item the level sets $W_s := \{z\in\Omega\ |\ \psi(z)< s\}$ are connected and relatively compact for all $s$ sufficiently negative.
\end{enumerate}
Additionally, we will require our weights to satisfy a differential inequality.
\begin{itemize}
    \item We say $\psi$ is an {\bf $m$-subweight} if $\psi$ is $m$-sh.
    \item We say $\psi$ is an {\bf $m$-superweight} if $\psi\in C^\infty(\Omega\setminus\Sing(\psi))$ and $(\ddbar \psi)^m\wedge\omega^{n-m} \leq 0$ there.
    \item We say $\psi$ is a {\bf maximal $m$-weight} if $\psi$ is $m$-sh such that $(\ddbar\psi)^m\wedge\omega^{n-m} = 0$ on $\Omega\setminus\Sing(\psi)$.
\end{itemize}

\end{definition}

Sub- and super-weights are suited to different measurements of the singularities of $m$-sh functions; subweights work well for Lelong numbers, while superweights are required when working with the relative type. Maximal $m$-weights are precisely those weights for which these two types of measurements can be compared.

\medskip

Now, by \cite{Dem_CADG_book, WanWang16}, for any $m$-subweight $\psi$, the complex Hessian operator:
\[
(\ddbar \psi)^m\wedge\omega^{n-m}
\]
is a well-defined, positive Borel measure on $\Omega$, with locally finite mass, which is continuous along decreasing sequences (that is to say that $\psi \in \mathcal{E}^m(\Omega)$, the $m$-sh version of the Cegrell class \cite{Lu15}). More generally, if $\psi_1, \ldots, \psi_{m-1}$ are all $m$-subweights and $\vp\in\mSH(\Omega)$ is arbitrary, then the mixed measure:
\[
\ddbar \vp \wedge \ddbar \psi_1\wedge \ldots \wedge \ddbar\psi_{m-1} \wedge\omega^{n-m}
\]
shares these same properties, and is moreover multilinear in each argument.

\begin{definition}\label{Integral_Lelong}
\cite{Dem_CADG_book, WanWang16} Suppose that $\psi$ is an $m$-subweight on $\Omega$ such that $\{\psi < -1\}\Subset \Omega$. Let $\vp\in \mSH(\Omega)$ be such that $\vp \leq -1$. For any $s < -1$, set:
\[
W_s := \{\psi < s\}.
\]

Then we define:
\begin{align}\label{eq_Integral_Lelong}
\nu_m(\vp, \psi) &:= \lim_{s\rightarrow -\infty} \int_{W_s} \ddbar \vp \wedge (\ddbar \psi)^{m-1}\wedge \omega^{n-m} \\
&= \int_{\{\psi = -\infty\}} \ddbar \vp\wedge(\ddbar \psi)^{m-1}\wedge \omega^{n-m}.\notag
\end{align}

We call $\nu_m(\vp, \psi)$ the {\bf generalized $m$-Lelong number} of $\vp$ with respect to $\psi$, or more compactly, the $(m, \psi)$-Lelong number of $\vp$. By standard results, $\nu_m(\vp, \psi)$ is always finite and non-negative, and the sequence in \eqref{eq_Integral_Lelong} is monotone decreasing.
\end{definition}

The following comparison theorem for $m$-polar measures is standard, and can be obtained by following the proof of \cite[Lemma 4.1]{ACCH09} (see also \cite{Dem93, HP17}):
\begin{proposition}\label{weight_compare_1}
Suppose that $\psi_1, \psi_2$ are $m$-subweights and $\vp_1, \vp_2\in \mSH(\Omega)$ are such that $\psi_1 \leq \psi_2 + C$ and $\vp_1\leq \vp_2 + C$ for some constant $C$. Let $S := \Sing(\psi_2)\cap\Sing(\vp_2)$. Then:
\begin{equation*}
\chi_{S} \ddbar \vp_2\wedge(\ddbar \psi_2)^{m-1}\wedge\omega^{n-m} \leq \chi_{S} \ddbar\vp_1\wedge(\ddbar \psi_1)^{m-1}\wedge\omega^{n-m}.
\end{equation*}
\end{proposition}

\medskip


\subsection{Relative Types}


Relative-types were first studied by Rashkovskii in \cite{Rash06}. The natural generalization of his original definition is the following lower semicontinuous function defined on $\Sing(\psi)$:

\begin{definition}\label{L_function}
Suppose that $\psi \leq -1$ is an $m$-superweight (or a maximal $m$-weight) and $\vp\in\mSH(\Omega)$. We define:
\[
L_\psi(\vp)(z_0) := \liminf_{\substack{z\rightarrow z_0 \\ z\not\in \Sing(\psi)}} \frac{\vp(z)}{\psi(z)}.
\]
It is easy to see that $L_{\psi}(\vp)$ is a lower-semincontinous function on $\Sing(\psi)$ (Proposition \ref{L_facts}).
\end{definition}

For a general $\psi$, this definition is somewhat lacking on its own. One of the main applications of the relative type is to get bounds of the form $\vp \leq L_{\psi}(\vp)\cdot \psi + C$ on a neighbourhood of $\Sing(\psi)$ -- but since $L_{\psi}(\vp)$ is non-constant, the right-hand side is only defined up to fixing an extension of $L_{\psi}(\vp)$ to a neighborhood of $\Sing(\psi)$. Since $L_{\psi}(\vp)$ is only lower semicontinuous. and $\Sing(\psi)$ can be very poorly behaved, it seems possible to us that this inequality may fail to hold if $L_{\psi}(\vp)$ is extended haphazardly.

A definition which can be more easily applied to get upper-bounds for $\vp$ is the following, which recovers the minimum of $L_{\psi}(\vp)$ (although we will see later that this definition has drawbacks of its own):

\begin{definition}\label{Relative_Type}
Suppose that $(\Omega^n, \omega)\subset (X^n, \omega)$ is an $m$-hyperconvex manifold. Let $\psi$ be an $m$-superweight (or maximal $m$-weight) on $\Omega$. Suppose that $W_0 = \{z\in\Omega\ |\ \psi(z)< 0\} \Subset \Omega$.

There exists a constant $A > 0$ such that for any $\vp\in\mSH(X,\omega)$, we have
\[
\vp+A\rho\in\mSH(\Omega).
\]
Define
\[
M_s(\vp+A\rho) := \max_{W_s}(\vp+A\rho).
\]
Then we define the {\bf relative type} of $\vp$ with respect to $\psi$ to be:
\begin{equation}\label{poorly_written_limit}
\sigma(\vp, \psi) := \lim_{s\rightarrow-\infty} \frac{M_s(\vp+A\rho)}{s}
= \lim_{s\rightarrow-\infty} \frac{\max_{W_s}(\vp)}{s}.
\end{equation}
\end{definition}
That $\sigma(\vp, \psi)$ is well-defined comes from the following three-circles type result (see e.g. \cite{Liu16} for similar results in the psh case):

\begin{proposition}\label{maximum_principle}
The function $M_s(\vp+A\rho)$ in Definition \ref{Relative_Type} is convex on $(-\infty, -1)$. It follows that the limit \eqref{poorly_written_limit} exists, and $\sigma(\vp, \psi)$ is always finite and non-negative.
\end{proposition}
\begin{proof}
Suppose first that $\psi$ is a superweight, and that $\vp$ is also smooth on $\Omega$. Let $\rho\leq 0$ be a strictly $m$-sh function on $\Omega$, and define $\vp_\e := \vp + (A + \e)\rho$. Let $-\infty < s_1 < s_2 < -1$, and consider the function:
\[
F(z) := \frac{s_2 - \psi(z)}{s_2 - s_1}M_{s_1}(\vp_\e) + \frac{\psi(z) - s_1}{s_2 - s_1}M_{s_2}(\vp_\e)
\]
defined on $W := W_{s_2}\setminus \ov{W}_{s_1}$. We seek to show that:
\[
\vp_\e(z) \leq F(z)\text{ on } W.
\]
It is clear this inequality holds on $\p W$, so suppose for the sake of a contradiction that $\vp_\e - F$ admits an interior maximum at some $z_0\in W$. It follows that:
\[
0 < (\ddbar \vp_\e(z_0))^m\wedge\omega^{n-m} \leq (\ddbar F(z_0))^m\wedge\omega^{n-m}.
\]
But this is impossible, as $F$ is a superweight. Taking the limit as $\e\rightarrow 0$ finishes the proof in this case.

The case of a general $\vp$ follows now by using a decreasing sequence of smooth $m$-sh functions which converge to $\vp$ \cite{LN15, Plis13}.

If $\psi$ is instead a maximal subweight, we may by-pass the smooth approximation argument and instead appeal to the comparison principle directly, since in this case, $F$ will be a maximal $m$-sh function on $W_{s_2}\setminus\ov{W}_{s_1}$.

\end{proof}

An immediate consequence of Proposition \ref{maximum_principle} is the following alternative definition for $\sigma(\vp, \psi)$:

\begin{proposition}\label{good_def}
$\sigma(\vp, \psi) = \max\{\gamma \geq 0\ |\ \vp \leq \gamma\,\psi + O(1)\}$.
\end{proposition}
\begin{proof}
The proof is standard -- it suffices to check that:
\[
\sigma(\vp, \psi) \in \{\gamma \geq 0\ |\ \vp\leq \gamma\,\psi + O(1)\}.
\]
Let $A > 0$ be such that $\vp + A\rho\in\mSH(\Omega)$. Convexity of $M_s(\vp + A\rho)$ implies that $\sigma(\vp, \psi)$ can be computed by the slopes of the secant lines:
\begin{equation}\label{limit}
\sigma(\vp, \psi) = \lim_{s\rightarrow-\infty}\frac{M_{s_0}(\vp + A\rho) - M_s(\vp + A\rho)}{s_0 - s}
\end{equation}
for any fixed $s_0 < -1$. Since $M_s(\vp + A\rho)$ is also decreasing as $s\rightarrow-\infty$, the sequence in \eqref{limit} is decreasing, and we have:
\[
M_s(\vp + A\rho) \leq \sigma(\vp, \psi) s + M_{s_0}(\vp  +A\rho) - \sigma(\vp, \psi) s_0,
\]
for all $s < s_0$. From this, we see that:
\[
\vp \leq \sigma(\vp, \psi) \psi + O(1)\quad \text{ on }W_{s_0}.
\]
\end{proof}

For later use, we record the following facts:
\begin{proposition}\label{L_facts}
The function $L_{\psi}(\vp)$ is lower semicontinuous on $S_\psi := \Sing(\psi)$, and:
\[
\sigma(\vp, \psi) = \min_{\Sing(\psi)} L_\psi(\vp).
\]
\end{proposition}
\begin{proof}
We first show the lower semicontinuity of $L_{\psi}(\vp)$. Fix a point $z_{0}\in S_\psi$ and let $z_{i}\in S_\psi$ be a sequence of points converging to $z_{0}$. For any $\e>0$, there exists $x_{i}\in X\setminus S_\psi$ such that
\[
\mathrm{dist}_{\omega}(x_{i},z_{i}) \leq i^{-1}, \quad
\frac{\vp(x_{i})}{\psi(x_{i})} \leq L_\psi(\vp)(z_{i})+\e.
\]
It is clear that $x_{i}\to z_{0}$. By definition:
\[
L_\psi(\vp)(z_{0}) \leq \liminf_{i\to\infty}\frac{\vp(x_{i})}{\psi(x_{i})} \leq \liminf_{i\to\infty}L_{\psi}(\vp)(z_{i})+\e.
\]
Letting $\e\rightarrow0$ concludes.

Now, it is clear from the definitions that $L_{\psi}(\vp) \geq \sigma(\vp, \psi)$. To see that it is the minimum, using that $\vp$ is upper semicontinuous, we can find a sequence of points $z_i\in W_{-i} = \{\psi \leq -i\}$ such that $\vp(z_i) = \max_{W_{-i}} \vp$. It follows from compactness of $S_\psi$ that there exists a convergent subsequence $z_i \rightarrow z_0\in S_\psi$, and so we see that $L_{\psi}(\vp)(z_0) \leq \sigma(\vp,\psi)$, finishing the proof.

\end{proof}

\medskip

Again following Demailly \cite{Dem_CADG_book}, Proposition \ref{weight_compare_1} can be used to compare $\sigma(\vp, \psi)$ and $\nu_m(\vp, \psi)$. For notational convenience, we set
\[
S_\psi := \Sing(\psi)\quad \text{ and }\quad\mu_\psi := \chi_{S_\psi}(\ddbar \psi)^m\wedge\omega^{n-m}.
\]
Note that $\Supp(\mu_\psi)\subseteq S_\psi$, but the inclusion may be proper; indeed, as Example 2.1 of \cite{ACH15} shows, $\Supp(\mu_\psi)$ can be quite small inside $S_\psi$.

\begin{corollary}\label{Lelong_Relative_Compare}
Suppose that $\psi \leq -1$ is a maximal $m$-weight such that $\mu_\psi(S_\psi) > 0$. Then:
\[
\sigma(\vp, \psi) \leq \frac{1}{\mu_\psi(S_\psi)}\nu_m(\vp, \psi).
\]
\end{corollary}
\begin{proof}
By Corollary \ref{good_def}, we have that:
\[
\vp \leq \sigma(\vp, \psi)\psi + C
\]
in some neighborhood of $S_\psi$. By Proposition \ref{weight_compare_1}, it follows that:
\[
\sigma(\vp, \psi) (\ddbar\psi)^m\wedge\omega^{n-m} \leq \ddbar\vp\wedge(\ddbar\psi)^{m-1}\wedge\omega^{n-m}.
\]
Integrating over $S_\psi$ finishes the proof.
\end{proof}

\begin{remark}
Proposition \ref{weight_compare_1} allows us to define the Radon-Nikodym derivative of $\chi_{S_\psi} \ddbar \vp\wedge(\ddbar\psi)^{m-1}\wedge\omega^{n-m}$ with respect to $\mu_\psi$ in certain cases, e.g. when $\vp\in \mathcal{E}^m(\Omega)$ (see \cite[Lem. 4.4]{ACCH09}). It is interesting to ask how different this function is from $L_{\psi}(\vp)$.
\end{remark}


\section{Construction of $m$-Weights Along Submanifolds}\label{section_constructions}


In this section, we construct our maximal $m$-weights $\psi_V$, associated to the complex submanifold $V$.


\subsection{Construction of the $m$-subweight associated to $V$}


Let $V\subseteq X$ be a smooth submanifold of codimension $k$. Let $r(z) := \mathrm{dist}_\omega(z, V)$ be the Riemannian distance function to $V$, and choose $0 < s_V < 1$ sufficiently small so that $\Omega := \{z\in X\ |\ r(z) < s_V\}$ does not intersect the cut-locus of $V$.

It follows that $r$ will be smooth on $\Omega\setminus V$, so that:
\[
r_\e := \sqrt{r^2 + \e}
\]
will be smooth on all of $\Omega$ for any $\e > 0$.

Define $G_m: \R_{\geq 0}\rightarrow \R\cup \{-\infty\}$:
\[
G_m(t) := \begin{cases} \log t &\text{ if } m = k \\ -t^{2 - \frac{2k}{m}} &\text{ if }m < k.\end{cases}
\]
$G_m$ will be the natural scaling of a maximal $m$-weight near $V$ -- by perturbing it slightly, we will be able to produce sub- and super-weights. Our perturbations will be of the form:
\[
h(s) := s + A s^{1 + \delta}.
\]
where $\max\left\{1, 2\left(\frac{k}{m} - 1\right)\right\} < \delta$ will be a uniform constant for all $s$ sufficiently small. We will also have $A = \pm 1$, depending on if we are constructing the sub- or super-solution. Note that as long as $\max\left\{1, 2\left(\frac{k}{m} - 1\right)\right\} < \delta$, $G_m(h(r))$ and $G_m(r)$ will have the same singularity type on $\Omega$, for either choice of $A$.\\

We will need to choose suitable coordinates for our computations, which we fix once and for all. Let $z\in \Omega\setminus V$, and let $x\in V$ be the unique nearest point to $z$. Consider the geodesic connecting these points, which will be normal to $V$ at $x$. Choose a unitary frame $\{e_{i}\}_{i=1}^{n}$ at $x$ such that
\[
e_{1} = \frac{1}{\sqrt{2}}(\nabla r - \sqrt{-1} J\nabla r),
\]
the $e_{2},\ldots,e_{k}$ are perpendicular to $V$, and the $e_{k+1},\ldots,e_{n}$ are tangent to $V$; we also write $\{e_{i}\}_{i=1}^{n}$ for the parallel transport of this frame to $z$. Then, according to \cite[Lemma 2.2]{TY12}, in these coordinates we have that:
\[
r_{i\ov{j}}=
\begin{pmatrix}
\frac{1}{2r} & 0 & 0 \\
0 & \frac{1}{r}\cdot\mathrm{Id}_{k-1} & 0 \\
0 & 0 & 0
\end{pmatrix}
+o(1),
\]
where $o(1)$ denotes a term satisfying $\lim_{r\rightarrow0}o(1)=0$ (see in particular equation (2.22) in \cite{TY12}).

In what follows, it suffices to restrict to the case $m \leq k$, since, if $m > k$, then the $m$-sh functions we construct will also be $k$-sh.

\begin{proposition}\label{full_expression}
Suppose that $z\in \Omega\setminus V$ and that $m \leq k$. Let $D_{k, m} = \left(\frac{2k}{m} - 2\right)^{-1}$ if $k\not= m$ and $D_{k,m} = 1$ if $k = m$.

Then for any $\e > 0$, we have:
\begin{align}
D_{k,m}h^{\frac{2k}{m}}(\ddbar G_m(h(r_\e)))_{i\ov{j}} &= \begin{pmatrix}
1 - \frac{k}{m} + AB_1r_\e^\delta & 0 & 0 \\
0 & \left(1 + AB_2 r_\e^\delta\right)\mathrm{Id}_{k-1} & 0 \\
0 & 0 & 0
\end{pmatrix}\notag \\
&+ \frac{\e}{r_\e^2}\begin{pmatrix} \frac{k}{m} + A B_3 r_\e^\delta + o(r_\e^\delta) & 0 & 0\\ 0& 0 & 0 \\ 0 & 0 & 0\end{pmatrix}
+ o(r_\e^\delta),\label{ugly}
\end{align}
where we define the constants:
\[
B_{1} =2\left[1-\frac{k}{m}+\left(1-\frac{k}{m}\right)\delta+\frac{\delta^{2}}{4}\right]
\]
and
\[
B_2 := 2 + \delta,\quad \quad B_3 := 2\left[\frac{k}{m} + \delta\left(\frac{k}{m} - \frac{1}{2}\right) - \frac{\delta^2}{4}\right].
\]
The terms $o(r_\e^\delta)$ satisfy $o(r_\e^\delta) \leq c r_\e^\delta$ for any $0 < c$ and all $0 < r_\e^\delta < r_c$, where $r_c$ depends on both $c$ and an upper bound for $\delta > \max\left\{1, 2\left(\frac{k}{m} - 1\right)\right\}$.
\end{proposition}

Since the constant $D_{k,m} > 0$, it will have no bearing on the positivity/negativity of $\ddbar G_m(h(r_\e))$, and we will omit it in all of the below computations.

For later use, we also remark that it will be clear from the below proof that \eqref{ugly} still holds if either $\e > 0$ and $z\in V$ or if $\e = 0$ and $z\not\in V$ (replacing $o(r_\e)$ with $o(r)$ in this case).

\begin{proof}
The proof is a computation. To start, if $m < k$, then we have:
\[
(G_m(h(r_\e)))_{i\ov{j}} = \left(\frac{2k}{m}-2\right)\cdot h^{-\frac{2k}{m}}\cdot\left( h h' (r_\e)_{i\ov{j}}  + \left(h h'' + \left(1-\frac{2k}{m}\right)(h')^2\right)(r_\e)_{i}(r_\e)_{\ov{j}}\right),
\]
while if $m = k$:
\[
(G_m(h(r)))_{i\ov{j}} = h^{-2}\cdot\left( h h' (r_\e)_{i\ov{j}}  + \left(h h'' -(h')^2\right)(r_\e)_{i}(r_\e)_{\ov{j}}\right).
\]
Thus, we may deal with both cases simultaneously by computing:
\begin{equation}\label{bounding}
 h h' (r_\e)_{i\ov{j}}  + \left(h h'' + \left(1-\frac{2k}{m}\right)(h')^2\right)(r_\e)_{i}(r_\e)_{\ov{j}}.
\end{equation}

We now compute:
\[
(r_{\e})_{i} = \frac{rr_{i}}{\sqrt{r^{2}+\e}} = \frac{rr_{i}}{r_{\e}}, \quad
(r_{\e})_{i}(r_{\e})_{\ov{j}} = \frac{r^{2}r_{i}r_{\ov{j}}}{r_{\e}^{2}}
\]
and
\[
\begin{split}
(r_{\e})_{i\ov{j}}
= {} & \frac{rr_{i\ov{j}}}{\sqrt{r^{2}+\e}}+\frac{r_{i}r_{\ov{j}}}{\sqrt{r^{2}+\e}}-\frac{r^{2}r_{i}r_{\ov{j}}}{(r^{2}+\e)^{\frac{3}{2}}}\\
={} & \frac{r r_{i\ov{j}}}{r_\e} + \frac{\e r_i r_{\ov{j}}}{r_\e^3}.
\end{split}
\]
Hence, in our coordinates, we have
\[
(r_{\e})_{i}(r_{\e})_{\ov{j}}
= \frac{r^{2}r_{i}r_{\ov{j}}}{r_{\e}^{2}} =
\begin{pmatrix}
\frac{r^{2}}{2r_{\e}^{2}} & 0 & 0 \\
0 & 0 & 0 \\
0 & 0 & 0
\end{pmatrix}
\]
and
\begin{align*}
(r_{\e})_{i\ov{j}} = \frac{rr_{i\ov{j}}}{r_{\e}}+\frac{\e r_{i}r_{\ov{j}}}{r_{\e}^{3}} &=
\begin{pmatrix}
\frac{1}{2r_{\e}}+\frac{\e}{2r_{\e}^{3}} & 0 & 0 \\
0 & \frac{1}{r_{\e}}\cdot\mathrm{Id}_{k-1} & 0 \\
0 & 0 & 0
\end{pmatrix}
+\frac{r}{r_\e}\cdot o(1).\\
\end{align*}
Applying this to \eqref{bounding} gives:
\begin{align*}
&  h h' (r_\e)_{i\ov{j}}  + \left(h h'' + \left(1-\frac{2k}{m}\right)(h')^2\right)(r_\e)_{i}(r_\e)_{\ov{j}} \\
= & \begin{pmatrix}
(h h')\left(\frac{1}{2r_\e}+\frac{\e}{2r_\e^3}\right) + \frac{r^2}{2r_\e^2}h h'' + \left(\frac{1}{2}-\frac{k}{m}\right)\frac{r^2}{r_\e^2}(h')^2 & 0 & 0 \\
0 & (h h')\frac{1}{r_\e}\cdot\mathrm{Id}_{k-1} & 0 \\
0 & 0 & 0
\end{pmatrix} + o(r)
\end{align*}
as $hh' = r_\e + o(r_\e)$.

We simplify this matrix. Recalling that $h(s) = s + As^{1+\delta}$, we have:
\[
h'(s) = 1 + A(1+\delta) s^{\delta},\ \ \text{ and }\ \
h''(s) = A\delta (1+\delta) s^{\delta-1}.
\]
It follows that:
\begin{align*}
h h' &= s + (2+\delta)A s^{1+\delta} + (1+\delta)A^2s^{1+2\delta}\\
h h'' &= (1+\delta)\delta A s^\delta  + (1+\delta)\delta A^2 s^{2\delta}\\
(h')^2 &= 1 + 2(1+\delta)As^\delta + (1+\delta)^2A^2 s^{2\delta},
\end{align*}
so that the first entry simplifies to
\begin{align*}
&(h h')\left(\frac{1}{2r_\e} + \frac{\e}{2r_\e^3}\right) + \frac{r^2}{2r_\e^2}h h'' + \left(\frac{1}{2}-\frac{k}{m}\right)\frac{r^2}{r_\e^2}(h')^2\\
&= \left(r_\e + (2+\delta)A r_\e^{1+\delta} + (1+\delta)A^2r_\e^{1+2\delta}\right)\left(\frac{1}{2r_\e} + \frac{\e}{2r_\e^3}\right)\\
&\quad \quad + \left(\frac{1}{2} - \frac{\e}{2r_\e^2}\right)\left((1+\delta)\delta A r_\e^\delta  + (1+\delta)\delta A^2 r_\e^{2\delta}\right)\\
&\quad \quad + \left(\frac{1}{2}-\frac{\e}{2r_\e^2}\right)\left(1 - \frac{2k}{m}\right)\left(1 + 2(1+\delta)Ar_\e^\delta + (1+\delta)^2A^2 r_\e^{2\delta}\right)\\
&= \left(\frac{1}{2} + \frac{\e}{2r_\e^2}\right)\left(1 + (2+\delta)A r_\e^{\delta} + (1+\delta)A^2r_\e^{2\delta}\right)\\
&\quad \quad + \left(\frac{1}{2} - \frac{\e}{2r_\e^2}\right)\left((1+\delta)\delta A r_\e^\delta  + (1+\delta)\delta A^2 r_\e^{2\delta}\right)\\
&\quad \quad + \left(\frac{1}{2} - \frac{\e}{2r_\e^2}\right)\left(1-\frac{2k}{m}\right)\left(1 + 2(1+\delta)Ar_\e^\delta + (1+\delta)^2A^2 r_\e^{2\delta}\right)\\
&= 1 - \frac{k}{m} + 2A\left[1 - \frac{k}{m} + \left(1 - \frac{k}{m}\right)\delta + \frac{\delta^2}{4}\right]r_\e^\delta + o(r_\e^{\delta})\\
&\quad \quad +\frac{\e}{r_\e^2} \left(\frac{k}{m} + 2A\left(\frac{k}{m} + \delta \left(\frac{k}{m} - \frac{1}{2}\right) - \frac{\delta^2}{4}\right)r_\e^\delta + o(r_\e^\delta)\right)\\
&= 1 - \frac{k}{m} + AB_1 r_\e^\delta + \frac{\e}{r_\e^2}\left(\frac{k}{m} + AB_3r_\e^\delta + o(r_\e^\delta)\right) + o(r_\e^\delta).
\end{align*}
As $\frac{h h'}{r_\e} = 1 + A B_2 r_\e^\delta + o(r_\e^\delta)$ and $0 < \max\{r, \sqrt{\e}\} \leq r_\e$, we conclude \eqref{ugly}.
\end{proof}


We may now construct our $m$-subweights.

\begin{proposition}\label{subweight_construction}
Suppose that $m\leq k$. Then for any $\delta > \max\left\{1, 2\left(\frac{k}{m} - 1\right)\right\}$, there exists a $0 < \ul{s} < s_V/2$ such that for all $0 < \e < s_V/2$ we have that:
\[
\underline{\psi}_{V, \e} := G_m(\ul{h}(r_\e)) \in \mSH(W_{\ul{s}}),
\]
where $\ul{h}(s) := s + s^{1 + \delta}$. Moreover, as $\e\rightarrow 0$, we have that
\[
\underline{\psi}_{V, \e} \searrow G_m(\ul{h}(r)) := \underline{\psi}_V \in \mSH(W_{\ul{s}}).
\]
The constant $\ul{s}$ depends on $k, m, \omega, V,$ and an upperbound for $\delta$.
\end{proposition}
\begin{proof}
Since $G_m(\ul{h}(r_\e))$ is smooth on $\Omega$, it will be sufficient to check that it is $m$-sh at an arbitrary $z\in \Omega\setminus V$ (though the computations translate easily to the case when $z\in V$). Choose coordinates as before, so that by Proposition \ref{full_expression} with $A = 1$ gives:
\begin{align}
D_{k,m}\ul{h}^{\frac{2k}{m}}(\ddbar \ul{\psi}_{V,\e})_{i\ov{j}} &= \begin{pmatrix}
1 - \frac{k}{m} + B_1r_\e^\delta & 0 & 0 \\
0 & \left(1 + B_2 r_\e^\delta\right)\mathrm{Id}_{k-1} & 0 \\
0 & 0 & 0
\end{pmatrix}\notag \\
&+ \frac{\e}{r_\e^2}\begin{pmatrix} \frac{k}{m} + B_3 r_\e^\delta + o(r_\e^\delta) & 0 & 0\\ 0& 0 & 0 \\ 0 & 0 & 0\end{pmatrix}
+ o(r_\e^\delta). \label{ugly_sub}
\end{align}
$B_3$ can be bounded using an upperbound on $\delta$, so as long as $\ul{s}$ is sufficiently small the second matrix will be positive semi-definite for any $\e > 0$, and can be dropped from all further computations.

We $k = m$, \eqref{ugly_sub} now simplifies to:
\begin{align*}
\ul{h}^{2}(\ddbar \ul{\psi}_{V,\e})_{i\ov{j}} &\geq \begin{pmatrix}
\frac{\delta^2}{2}r_\e^\delta & 0 & 0 \\
0 & \mathrm{Id}_{k-1} & 0 \\
0 & 0 & 0
\end{pmatrix} + o(r_\e^\delta).
\end{align*}
The error term can be controlled by shrinking $\ul{s}$ if necessary, again depending on $\delta$, completing this case.

We now deal with the case when $m < k$. We compute the leading term of the $j$-th symmetric polynomial ($1 \leq j \leq m$) of \eqref{ugly_sub} to be:
\begin{align*}
&\left(1 - \frac{k}{m} + B_1 r_\e^{\delta}\right)(1+ B_2r_\e^{\delta})^{j-1} + \left(\frac{k}{j}-1\right) (1+ B_2r_\e^{\delta})^j\\
&= \frac{k}{j} - \frac{k}{m} + \left[B_1 + \left( \left(1-\frac{k}{m}\right)(j-1) + \left(\frac{k}{j} - 1\right)j\right)B_2\right]r_\e^{\delta} + o(r_\e^{\delta}),
\end{align*}
omitting a multiplicative factor of $j \cdot (k-1)\ldots(k-j+1)$. When $j < m$, the leading term is positive. When $j = m$, this term is zero, so the leading order term becomes:
\begin{align*}
&\left[B_1 + \left( \left(1-\frac{k}{m}\right)(m-1) + \left(\frac{k}{m} - 1\right)m\right)B_2\right]r_\e^{\delta}\\[1mm]
= {} & \left[B_1 + \left(\frac{k}{m}-1\right)B_2\right]r_\e^{\delta}\\[1mm]
= {} & \delta\left[\frac{\delta}{2}-\left(\frac{k}{m}-1\right)\right]r_\e^{\delta},
\end{align*}
since we have assumed:
\[
\delta > 2\left(\frac{k}{m} - 1\right).
\]
By taking $\ul{s}$ sufficiently small depending on the above constant, we can control the error term, showing that $\ddbar\ul{\psi}_{V,\e}$ is an $m$-subweight. Letting $\e\rightarrow 0$ shows the secondary statement in the proposition immediately.
\end{proof}


\subsection{Construction of the $m$-superweight}


We may now construct our $m$-superweight by very similar considerations:
\begin{proposition}\label{superweight_construction}
Suppose that $m\leq k$. Then for any $\delta > \max\left\{1, 2\left(\frac{k}{m} - 1\right)\right\}$, there exists a $0 < \ov{s} < s_V/2$, depending on $k, m, \omega, V,$ and an upperbound for $\delta$, so that if $\ov{h}(s) := s - s^{1 + \delta}$, then:
\[
\ov{\psi}_V := G_m(\ov{h}(r))
\]
is an $m$-superweight on $W_{\ov{s}}$.
\end{proposition}
\begin{proof}
By the remark immediately following Proposition \ref{full_expression}, we have:
\begin{align}\label{ugly_sup}
\ov{h}^{\frac{2k}{m}}(\ddbar \ov{\psi}_V)_{i\ov{j}} &= \begin{pmatrix}
1 - \frac{k}{m} - B_1r^\delta & 0 & 0 \\
0 & \left(1 - B_2 r^\delta\right)\mathrm{Id}_{k-1} & 0 \\
0 & 0 & 0
\end{pmatrix} + o(r^\delta).
\end{align}
When $k = m$, \eqref{ugly_sup} simplifies to:
\begin{align*}
\ov{h}^2(\ddbar \ov{\psi}_V)_{i\ov{j}} &= \begin{pmatrix}
-\frac{\delta^2}{2}r^\delta & 0 & 0 \\
0 & (1 - (2 + \delta)r^\delta)\mathrm{Id}_{k-1} & 0 \\
0 & 0 & 0
\end{pmatrix} + o(r^\delta).
\end{align*}
Choosing $\ov{s}$ sufficiently small relative to $\delta$ again controls the error term, so that this is a superweight.

When $m < k$, we again compute the leading terms of the $j$-th symmetric polynomials of the matrix in \eqref{ugly_sup}. When $j < m$, the leading order term is the same as the leading term in Proposition \ref{subweight_construction}; when $m = j$, it becomes the negation of the term in Proposition \ref{subweight_construction}, and hence will be negative as long as $\delta > 2\left(\frac{k}{m} - 1\right)$. Again, the error terms are an order of magnitude smaller than the leading order term, and so can be controlled assuming an upper bound for $\delta$.
\end{proof}

\subsection{The maximal $m$-weight along $V$}


We now show Theorem \ref{psi_V_intro}, and construct $\psi_V$ from $\ul{\psi}_V$ and $\ov{\psi}_V$ by using an envelope to make the complex Hessian measure vanish.

\begin{theorem}\label{maximal_construction}
Suppose that $m\leq k$, and $\delta > \max\left\{1, 2\left(\frac{k}{m} - 1\right)\right\}$. Let $\ul{\psi}_V$ and $\ov{\psi}_V$ be the $m$-subweight and $m$-superweights constructed in Propositions \ref{subweight_construction} and \ref{superweight_construction}, respectively, for the given $\delta$. Set $s_0 := \min\{\ul{s}, \ov{s}\}$. Then the function:
\[
\psi_V := \sup\{\vp\in\mSH(X,\omega)\ |\ \vp \leq 0 , \vp(z)\leq G_m(r(z)) + O(1)\}^*,
\]
is a maximal $m$-weight on $W_{s_0}$, which moreover has the same singularity type as $\ul{\psi}_V, \ov{\psi}_V$, and $G_m(r)$.
\end{theorem}
\begin{proof}
As mentioned at the beginning of this section, the three functions $\ul{\psi}_V, \ov{\psi}_V,$ and $G_m(r)$ have the same singularity type, so that
\[
\ul{\psi}_V - C_0 \leq G_m(r) \leq \ov{\psi}_V + C_0 \leq \ul{\psi}_V + C_0
\]
for some sufficiently large constant $C_0$. Since it is clear that the definition of $\psi_V$ only depends on the singularity type of $G_m(r)$, we have that:
\[
\psi_V := \sup\{\vp\in\mSH(X,\omega)\ |\ \vp \leq 0 , \vp(z)\leq \ov{\psi}_V(z) + O(1)\}^*.
\]
Since $\ul{\psi}_V \leq 0$, the envelopes:
\[
\psi_{V, C} := \sup\{\vp\in\mSH(X,\omega)\ |\ \vp \leq \min\{0, \ov{\psi}_V + C\}\}^*
\]
satisfy:
\[
\ul{\psi}_V \leq \psi_{V, C}
\]
for all $C \geq C_0$. Recall that the $\psi_{V,C}$ increase to $\psi_V$ as $C\rightarrow \infty$, so that $\ul{\psi}_V \leq \psi_{V}$ as well.

Let $A > 0$ be such that $\ov{\psi}_V + A \geq 0$ on $\p W_{s_0}$. By Lemma \ref{ugasdf} below, we have that:
\[
\psi_{V, C} \leq \ov{\psi}_V + A \text{ on } W_{s_0}.
\]
Letting $C\rightarrow \infty$ now shows that $\psi_{V} \leq \ov{\psi}_V + A$, so $\psi_V$ has the same singularity type as $G_m(r)$.

The conclusion about vanishing mass follows from a standard balayage argument -- see e.g. \cite{BT76}.

\end{proof}

The following lemma is just a restatement of \cite[Lemma 4.1]{RS15}.

\begin{lemma}\label{ugasdf}
Suppose that $\psi\leq 0$ is an $m$-superweight (or a maximal $m$-weight) and $\vp\in\mSH(X, \omega)$, $\vp \leq 0$, is such that $\sigma(\vp, \psi) \geq 1$. If $\vp \leq \psi \text{ on }\p\Omega$, then there exists a constant, independent of $\vp$, so that $\vp \leq \psi + C$ on all of $\Omega$.
\end{lemma}
\begin{proof}
Let $s_0 \leq 0$ be such that $W_{s_0} := \{\psi(z) < s_0\} \Subset \Omega$. Since $\sigma(\vp, \psi) \geq 1$, we can use Proposition \ref{maximum_principle} to see that:
\[
(\vp+A\rho)(z) \leq \psi(z) - s_0 + M_{s_{0}}(\vp+A\rho) \leq \psi(z) - s_0,
\]
for all $z\in W_{s_0}$.

It follows that $\vp+A\rho\leq \psi - s_{0}$ on $\p(\Omega\setminus W_{s_0})$. We can now conclude by the maximum/comparison principle, as in Proposition \ref{maximum_principle} (depending on if $\psi$ is a super/maximal weight).
\end{proof}

We now show Corollary \ref{Lelong_number_intro}, proving a formula for $\nu(\vp, \psi_V)$ as a density of a weighted Laplacian of $\vp$, which has been averaged over $V$.

\begin{proposition}\label{integral_Lelong_Laplacian}
Suppose that we are in the setting of the previous section, so that $r$ is the geodesic distance to $V^{n-k}\subset \Omega^n$. Let $\vp\in\mSH(\Omega)$ with $m \leq k$. Then:
\[
\nu_{m}(\vp, \psi_V) = \lim_{s\rightarrow 0}\frac{C_{k,m}}{s^{2k - 2k/m}}\int_{\{r < s\}} \ddbar \vp \wedge (\ddbar r^2)^{m-1}\wedge\omega^{n-m}.
\]
for some constant $C_{k, m}$ depending only on $k, m$.
\end{proposition}
\begin{proof}
For the sake of convenience, we only give the proof when $k > m$. Let $H_m(t) := \left(-t\right)^{-\frac{m}{k-m}}$. Note that $H_m$ is convex and increasing on $(-\infty, 0)$, and that:
\[
H_m(G_m(t)) = t^2.
\]
Since $\{\ul{\psi}_V < G_m(s)\} = \{\ul{h}(r) < s\}$, by \cite[Chapter III, Formula (5.5)]{Dem_CADG_book}, we have that:
\begin{align*}
\int_{\{\ul{h} < s\}} &\ddbar\vp\wedge (\ddbar \ul{h}^2)^{m-1}\wedge\omega^{n-m}
\\&= (H'_m(G_m(s)))^{m-1} \int_{\{\ul{\psi}_V \leq G_m(s)\}} \ddbar \vp\wedge (\ddbar \ul{\psi}_V)^{m-1}\wedge \omega^{n-m}\\
&=\left(\frac{m}{k-m}\right)^{m-1}\left(s^{\frac{2(m-k)}{m}}\right)^{\frac{-k(m-1)}{k-m}}\int_{\{\ul{\psi}_V \leq G_m(s)\}} \ddbar \vp\wedge (\ddbar \ul{\psi}_V)^{m-1}\wedge \omega^{n-m},
\end{align*}
so that:
\[
\nu_{m}(\vp, \ul{\psi}_V) = \left(\frac{k-m}{m}\right)^{m-1}\lim_{s\rightarrow 0}\frac{1}{s^{2k-2k/m}}
\int_{\{\ul{h} < s\}} \ddbar\vp\wedge (\ddbar \ul{h}^2)^{m-1}\wedge\omega^{n-m}.
\]

Since $\ul{h}$ is only a small perturbation of $r$, we can replace it with $r$ in the above expression by the following argument. By \cite{TY12}, $\ddbar r^2 = \begin{pmatrix}I_{k} & 0 \\ 0 & 0 \end{pmatrix} + o(1)$ is smooth and $m$-sh on any sufficiently small neighborhood of $V$. By possibly increasing $\delta$, we may assume that it is a large integer, so that $\ul{h}^2$ is also smooth.

Let $T := \ddbar \vp \wedge (\ddbar \ul{h}^2)^{m-2}\wedge\omega^{n-m}$. We have that:
\[
\ddbar \ul{h}^2 = \ddbar \left(r + r^{1 +\delta}\right)^2 = \ddbar \left(r^2 + 2r^{2+\delta} + r^{2+2\delta}\right),
\]
so away from $V$:
\begin{align*}
\ddbar r^{2+\delta} & = (2+ \delta)\sqrt{-1}\p \left( r^{1+\delta} \pbar r\right)\\[1mm]
& = (2+\delta)r^{1+\delta}\ddbar r + (2 + \delta)(1+\delta)r^{\delta}\sqrt{-1}\p r\wedge \pbar r\\[1mm]
 & = \frac{2+\delta}{2} r^\delta \ddbar r^2 + \delta(2+\delta) r^\delta\sqrt{-1}\p r\wedge \pbar r.
\end{align*}
Since both sides are continuous across $V$ however, the above expression holds there as well. Then:
\begin{align*}
\ddbar \ul{h}^2 - \ddbar r^2
= {} & \left((2+\delta)r^\delta + (1+\delta)r^{2\delta}\right)\ddbar r^2  \\
& + 2\delta\left((2+\delta)r^\delta + (2+2\delta)r^{2\delta}\right)\sqrt{-1}\p r\wedge \pbar r
\end{align*}
and so
\[
0 \leq \ddbar \ul{h}^2 - \ddbar r^2 \leq C_{\delta}r^{\delta}\ddbar \ul{h}^2,
\]
where the above $\leq$'s are understood as meaning that the difference between the two sides is $m$-sh. It follows that:
\begin{align*}
0 &\leq\frac{\chi_{\{\ul{h} < s\}}}{s^{2k - 2k/m}}\left( T\wedge \ddbar \ul{h}^2 - T\wedge \ddbar r^2\right)\\
& \leq C_{\delta}s^\delta \left(\frac{\chi_{\{\ul{h} < s\}}}{s^{2k - 2k/m}} T\wedge \ddbar \ul{h}^2\right)
\end{align*}
since the measure $T\wedge \ddbar\ul{h}^2 \geq 0$. Taking $s\rightarrow 0$ now implies that:
\[
\lim_{s\rightarrow 0} \frac{1}{s^{2k - 2k/m}}\int_{\{\ul{h} < s\}} T \wedge \ddbar \ul{h}^{2}
= \lim_{s\rightarrow 0} \frac{1}{s^{2k - 2k/m}}\int_{\{\ul{h} < s\}} T \wedge \ddbar r^2,
\]
and it is easy to see that, by repeating the above argument $m-2$ more times, we have that:
\[
\nu(\vp, \ul{\psi}_V) = \left(\frac{k-m}{m}\right)^{m-1}\lim_{s\rightarrow 0} \frac{1}{s^{2k - 2k/m}}\int_{\{\ul{h} < s\}}\ddbar \vp \wedge (\ddbar r^2)^{m-1}\wedge \omega^{n-m}.
\]
By noting that $\{r < (1-\e)s\} \subset \{\ul{h} < s\} \subset \{r < s\}$ as long as $s$ is sufficiently small (depending on $\e$), we can also replace the set $\{\ul{h} < s\}$ with $\{r < s\}$ in the above limit.

Finally, since $\ul{\psi}_V$ and $\psi_V$ have the same singularity type, we can conclude by using Proposition \ref{weight_compare_1}.
\end{proof}


\subsection{Subweights with localized singularities}\label{subsection_localization}


We now construct explicit $m$-subweights, $\psi_\theta$, which have ``non-constant" behavior along $V$ when $k \geq m$. It is easy to see from the construction that these subweights will satisfy $\sigma(\psi_\theta, \psi_V) = 0$ and $\nu(\psi_\theta, \psi_V)  > 0$ (Remark \ref{Polar_Mass}).

For any real number $\nu \in \R$, define $F_\nu: \R_{\geq 0}\rightarrow \R\cup \{-\infty\}$ by:
\[
F_\nu(t) := \begin{cases} t^{2 - 2\nu} &\text{ if }\nu < 1 \\ \log t &\text{ if } \nu = 1 \\ -t^{2 - 2\nu} &\text{ if }\nu > 1.\end{cases}
\]
We additionally define $D_\nu := |2\nu - 2|$ if $\nu\not=1$, and $D_\nu = 1$ if $\nu = 1$. If $m \leq k$, then for any $\nu < \frac{k}{m}$ we have that $G_m < F_{\nu}$.

\begin{proposition}\label{localization}
Suppose we are in the setting of the previous subsections, and $\delta$ is sufficiently large, depending only on $k$ and $m$. Let $\theta$ be a smooth, non-negative function on $V$, and suppose that $m \leq k$. Then for any $\frac{k}{m} - \frac{1}{2} \leq \nu < \frac{k}{m}$, there exists a $C > 0$ such that:
\[
\ul{\psi}_{\theta} := \ul{\theta\psi}_V + C F_{\nu}(\ul{h}(r))\  \in\,\mSH(\Omega),
\]
where $\ul{\theta}$ is a smooth extension of $\theta$ with compact support. When $m = 1$, we can take $k - 1\leq\nu < k$.

The construction is local, and so works if $V$ is only locally defined in $\ov{\Omega}$.
\end{proposition}

\begin{remark}\label{ACCH_remark}
The bounds on $\nu$ seem to be optimal. Note that when $k < \frac{3}{2}m$, we can choose $\nu < 1$, so that the $F_{\nu}$ term is bounded from below.
\end{remark}

\begin{proof}
The proof is a computation which is very similar to Proposition \ref{full_expression}, so we omit several details. Working in the same geodesic normal coordinates as before, it follows that:
\begin{align*}
\ddbar(\ul{\theta\psi}_V(\ul{h}(r_\e)))
= {} & \ul{\theta}\ddbar G_m(\ul{h}(r_\e)) + \sqrt{-1} \p \ul{\theta} \wedge \pbar G_m(\ul{h}(r_\e))\\[3mm]
& +\sqrt{-1}  \p G_m(\ul{h}(r_\e))\wedge \pbar \ul{\theta} + G_m(\ul{h}(r_\e)) \ddbar\ul{\theta}\\[3mm]
\geq {} & \begin{pmatrix} 0 & 0 & \frac{D_{k,m}^{-1} r}{r_\e^{\frac{2k}{m}}}\pbar \ul{\theta} \\
0 & 0 & 0 \\
\frac{D_{k,m}^{-1} r}{r_\e^{\frac{2k}{m}}}\p \ul{\theta} & 0 & \frac{-1}{r_\e^{\frac{2k}{m}-2}}\ddbar \ul{\theta} \end{pmatrix},
\end{align*}
where we again understand $\geq$ as meaning the difference of the two sides is $m$-sh in this proof. Here are also interpreting $\pbar\ul{\theta}$ as a $1\times k$ matrix and $\p\ul{\theta}$ as its conjugate transpose.

From the previous proposition, we also have:
\begin{align*}
(\ddbar F_{\nu}(\ul{h}(r_\e)))_{i\ov{j}} &\geq D_{\nu}^{-1}\ul{h}^{-2\nu}\begin{pmatrix}
1 - \nu & 0 & 0 \\
0 & \mathrm{Id}_{k-1} & 0 \\
0 & 0 & 0
\end{pmatrix} + O(r_\e^{\delta - 2\nu}).
\end{align*}
This leading order matrix will be strictly $m$-sh as long as $\nu < \frac{k}{m},$ by previous computations.

We now have that
\[
\ddbar(\ul{\theta\psi}_V + C F_{\nu}(\ul{h}(r)))\geq \begin{pmatrix} \frac{C}{r_\e^{2\nu}}\left(1 - \nu\right) & 0 & \frac{r}{r_\e^{\frac{2k}{m}}}\pbar \ul{\theta} \\
0 & \frac{C}{r_\e^{2\nu}}\mathrm{Id}_{k-1} & 0 \\
\frac{r}{r_\e^{\frac{2k}{m}}}\p \ul{\theta} & 0 & \frac{-1}{r_\e^{\frac{2k}{m}-2}}\ddbar \ul{\theta} \end{pmatrix} + O(r_\e^{\delta - 2\nu}),
\]
ignoring some multiplicative constants. This leading matrix will be (strictly) $m$-sh with a positive term of order $r^{2\nu}_\e$ if we choose $\nu$ appropriately. This can be done by choosing $\nu$ so that the strict $m$-subharmonicity of the $F_\nu$ term beats the powers coming from the $\ul\theta$ terms. This can be seen by noting that $r < r_\e$, so that, when $m > 1$, we need:
\[
2\nu  \geq \frac{2k}{m} -1,
\]
giving the combined constraints:
\[
\frac{k}{m}  - \frac{1}{2} \leq \nu < \frac{k}{m};
\]
when $m = 1$, we only need the trace to be positive, so we only require:
\[
2\nu \geq 2k - 2 \implies k - 1\leq \nu < k.
\]

\end{proof}

\begin{remark}\label{Polar_Mass}
For $\vp$ in Cegrell's class $\mathcal{E}^m$, we denote the complex Hessian operator of $\vp$ by $\Hm(\vp) := \left(\ddbar\vp\right)^m\wedge\omega^{n-m}$. We claim that the $\ul{\psi}_\theta \in \mathcal{E}^m$ constructed in Proposition \ref{localization} satisfies:
\[
\chi_{V} \ddbar\ul{\psi}_{\theta} \wedge (\ddbar \psi_V)^{m-1}\wedge\omega^{n-m}= \chi_V\theta\Hm\left(\psi_V\right).
\]
Indeed, on any sufficiently small open $B$ near $z\in V$, we have:
\[
(\theta(z) + \e)\psi_V - C \leq \ul{\psi}_\theta \leq (\theta(z) - \e)\psi_V + C,
\]
so Proposition \ref{weight_compare_1} implies that
\[
\chi_{V\cap B}(\theta(z) - \e)^m\Hm\left(\psi_V\right)\leq \chi_{V\cap B} \Hm(\ul{\psi}_{\theta}) \leq \chi_{V\cap B}(\theta(z) + \e)^m\Hm\left(\psi_V\right).
\]
Covering $V$ with a disjoint collection of such $B$ and then letting $\e\rightarrow 0$ produces the claimed equality. 

Taking $\theta$ such that $\min_V \theta = 0$ but $\theta \not= 0$ now gives a function such that:
\[
\sigma(\ul{\psi}_\theta, \psi_V) = 0,\quad \text{ but }\quad \nu(\ul{\psi}_\theta, \psi_V) > 0. 
\]

\noindent Additionally, we can use Proposition \ref{weight_compare_1} to see that $\chi_{V} \Hm(\ul{\psi}_{\theta}) = \chi_V\theta^m\Hm\left(\psi_V\right).$ Solutions to $\Hm(v) = \chi_V \theta^m \Hm(\vp)$ for any given $\vp\in \mathcal{E}^m(\Omega)$ were constructed in \cite{HP17}, following the work of \cite{ACCH09} in the psh case, as a nested sequence of envelopes. Our $\ul{\psi}_\theta$ can be seen as explicit contenders in these envelopes for $\vp = \psi_V$.

Let $\psi_\theta$ be the solution to $\Hm(\psi_\theta) = \chi_V\theta^m\Hm(\psi_V)$ constructed in \cite{HP17}. If $k \leq \frac{3}{2}m$, our examples show that $\psi_\theta$ is actually bounded on $V\setminus \mathrm{Supp}(\theta)$. We suspect our examples actually compute the exact singularity type of $\psi_\theta$ on all of $V$.
\end{remark}


\subsection{Subharmonic Functions and Minimal Submanifolds}\label{subsection_minimal}


We conclude this section with the observation that the above method can be copied directly when $V$ is only a minimal submanifold of $X$ (which we still assume to be K\"ahler), provided that we take $m=1$. Subharmonic functions have long been known to be intimately related with minimal submanifolds, and as such, we are unsure if a similar construction is already present in the literature.

\begin{proposition}
Suppose that $V$ is a minimal submanifold of $X$ of (real) codimension $k$. Then there exists an $\omega$-subharmonic function $\psi_V$ on a neighborhood $\Omega$ of $V$, which is moreover $\omega$-harmonic on $\Omega\setminus V$, such that $\psi_V$ has the same singularity type as:
\[
-\frac{1}{r^{k - 2}}.
\]
\end{proposition}
\begin{proof}
Construct real coordinates at a point $z$ near $V$ in the same way as done at the beginning of this section, so that $e_1$ is in the direction of the geodesic connecting $z$ and $V$. By \cite[Lemma 2.2]{TY12}, we have:
\[
r_{ij}=
\begin{pmatrix}
0 & 0 & 0 \\
0 & \frac{1}{r}\cdot\mathrm{Id}_{k-1} & 0 \\
0 & 0 & h_{ij}
\end{pmatrix}
+o(1),
\]
where here $(h_{ij})_{k+1\leq i,j \leq 2n}$ are the components of the second fundamental form of $V$ evaluated against $e_1$. All of the computations in Proposition \ref{full_expression} go through verbatim, with the exception that the second fundamental form term becomes $rh_{ij} + o(r_\e)$. Since $V$ is minimal, this bad term disappears when we take the trace of $\ddbar r_\e^{k-2}$, and the computations in Propositions \ref{subweight_construction} and \ref{superweight_construction} can then be repeated.

Once a sub- and super-weight have been produced, Proposition \ref{maximal_construction} is standard.
\end{proof}


\section{A Siu-Type Theorem}\label{section_Siu}


We now present our Siu-type theorem, stated as Theorem \ref{Siu_intro} in the introduction. We will need the following weak Harnack inequality:

\begin{proposition}[Theorem 8.18 of \cite{GT83}]\label{weak_Harnack}
Let $\Omega\subset\mathbb{R}^{2n}$ be a domain. Suppose that $\vp$ is a non-positive smooth subsolution of a linear elliptic operator $L = a^{ij}\de_{i}\de_{j}$, i.e.
\[
\vp \leq 0, \quad L\vp \geq 0 \quad \text{in $\Omega$},
\]
and $\lambda I_{2n} \leq (a^{ij}) \leq \Lambda I_{2n}$ for positive constants $\lambda$ and $\Lambda$. Then for any ball $B_{2}(y)\subset\Omega$ and $p\in[1,n/(n-1))$, we have
\[
\vp(y) \leq -\frac{1}{C(\lambda,\Lambda,p,n)}\int_{B_{1}(y)}(-\vp)^{p} \leq 0.
\]
\end{proposition}

\medskip

\begin{theorem}\label{Siu}
Suppose that $V$ is a compact, complex submanifold of codimension $k$. If $\vp\in\mSH(\Omega)$ is such that $k < m$, then $L_{\psi_V}(\vp) \equiv \sigma(\vp, \psi_V)$.
\end{theorem}
\begin{proof}
For convenience, set $\sigma := \sigma(\vp, \psi_V)$. We will show that the set of points where $L_{\psi_V}(\vp) = \sigma$ is open in $V$ -- since $L_{\psi_V}(\vp)$ is lower semi-continuous, it is clear that this set is closed and non-empty, by Proposition \ref{L_facts}.

Let $z_0\in V$ be such that $L_{\psi_V}(\vp)(z_0) = \sigma$. Choose local holomorphic coordinates $z = (z', z'')$ centered at $z_0$ such that $V = \{z' = 0\}$ in these coordinates. Let $B_{R}$ be the coordinate ball of radius $0 < R\leq 2$, centered at $z_0 = (0, 0)$. Fix a uniform constant $C$ such that:
\begin{equation}\label{equivalence of r}
C^{-1} r(z) \leq |z'| \leq C r(z)\quad\text{ for all }z\in B_2.
\end{equation}
Further, let $\rho(z)$ be a smooth potential for $\omega$ on $B_2$ such that $-C \leq \rho(z) \leq 0$.

By definition, there exists a sequence $z_i = (z_i', z_i'') \rightarrow z_0$, with $z_i \in \Omega\setminus V$, such that:
\[
\lim_{i\rightarrow\infty} \frac{\vp(z_i)}{\psi_V(z_i)} = \sigma.
\]
We assume without loss of generality that $\log|z_i'|, \log r(z_i), \psi_V(z_i) \leq -1$. Since $\psi_V$ and $\log r$ have the same singularity type, by \eqref{equivalence of r}, we have that:
\[
\lim_{i\rightarrow\infty} \frac{\vp(z_i)}{\log |z_i'|} = \sigma
\]
also. It follows that, for any $\e > 0$, we have:
\begin{equation}\label{bounds for vp}
(\sigma + \e) \log |z_i'| \leq \vp(z_i) \leq (\sigma - \e)\log |z_i'|
\end{equation}
for all $i$ sufficiently large.

Define the functions $\eta_{i, \e}$ on $V_2 := V\cap B_2$ by:
\[
\eta_{i, \e}(z'') := \frac{\vp(z_i', z'') + \rho(z_i', z'') - \vp(z_i) + 3\e\log|z_i'|}{-\log|z_i'|}.
\]
Note that $\eta_{i, \e}$ is subharmonic on $V_2$ with respect to the metric $\omega|_{\{z' = z'_i\}}$ since we have assumed that $m > k$.

We will show that $\eta_{i_j, \e_j}\rightarrow 0$ pointwise a.e. for some subsequence $i_j\rightarrow\infty$ of the $i$ and an associated sequence of $\e_j\rightarrow0$, by using the Harnack inequality. To utilize this, we need to first verify that $\sup \eta_{i, \e} \leq 0$. Increasing $C$ if necessary, we have that:
\[
\{z' = z_i'\} \subseteq \{\log |z'| \leq \log |z_i'|\} \subseteq  \{\psi_V(z) \leq \log |z_i'| + C \},
\]
so by the definition of $\sigma$ we have:
\[
\sup_{V_2} \vp(z_{i}',\cdot) \leq \sup_{\{\psi_V(z) \leq \log |z_i'| + C\}} \vp \leq (\sigma - \e)(\log |z_i'| + C).
\]
Combining this with \eqref{bounds for vp} and our choice of $\rho$ gives:
\begin{equation}\label{eta i e negative}
\sup_{V_2}\eta_{i,\e} \leq -(\sigma - \e) + (\sigma + \e) - 3\e + \frac{(\sigma - \e)C}{-\log|z_i|} \leq -\e + \frac{(\sigma - \e)C}{-\log|z_i|} < 0
\end{equation}
for $\e$ fixed and $i$ sufficiently large.

Let $\e_{j}>0$ be a sequence such that $\e_{j}\to0$ as $j\to\infty$. By \eqref{eta i e negative}, for each $j$ there exists a sufficiently large $i_{j}$ such that $i_{j}\to\infty$ as $j\to\infty$ and
\[
\sup_{V_2}\eta_{i_{j},\e_{j}}<0.
\]
Letting $V_1 := V\cap B_1$, Proposition \ref{weak_Harnack} now implies that:
\[
0 \leq \int_{V_1} (-\eta_{i_{j}, \e_{j}})^p \leq -C \eta_{i_{j}, \e_{j}}(z_{i_{j}}'') \leq C\e_{j},
\]
for any $p \in [1, n/(n-1))$ (up to increasing $C$ again, depending on $p$ and $\omega$). Fixing some $p>1$, it is clear that
\[
\lim_{i\to\infty}\norm{\eta_{i_{j},\e_{j}}}_{L^p(V_1)} = 0.
\]
Hence, we see that $\eta_{i_{j},\e_{j}}\rightarrow 0$ in $L^p$ as $j\rightarrow\infty$, and so after passing to a subsequence of the $j$, we have that:
\begin{equation}\label{eta_to_zero}
\eta_{i_{j},\e_{j}}\rightarrow 0 \ \ \text{pointwise a.e. on $V_1$.}
\end{equation}

\medskip

Now, from our definition of $\eta_{i, \e}$, we see that:
\[
\frac{\vp(z'_{i_{j}}, z'')}{\log|z_{i_{j}}'|} = \frac{\vp(z_{i_{j}}', z_{i_{j}}'')}{\log|z_{i_{j}}'|}
-\eta_{{i_{j}}, \e_{j}}(z'') - 3\e_{j} - \frac{\rho(z'_{i_{j}}, z'')}{\log|z'_{i_{j}}|}.
\]
Using \eqref{eta_to_zero} and taking the limit as $j\rightarrow\infty$ gives:
\[
\lim_{j\rightarrow\infty} \frac{\vp(z'_{i_{j}}, z'')}{\log|z_{i_{j}}'|} = \sigma
\]
for almost all $z''\in V_1$. Then we conclude that:
\[
L_{\psi_V}(\vp)(0, z'') \leq \lim_{j\rightarrow\infty} \frac{\vp(z'_{i_{j}}, z'')}{\log|z_{i_{j}}'|} = \sigma
\]
for almost every $z''\in V_1$. It follows now from Proposition \ref{L_facts} that in fact $L_{\psi_V}(\vp) \equiv \sigma$ on $V_1$, finishing the proof.
\end{proof}


\end{document}